\documentclass[a4paper,12pt]{amsart}
\usepackage{amsfonts}
\usepackage{amssymb}
\usepackage{ifthen}
\usepackage{amscd}
\usepackage{amsxtra}
\usepackage{graphicx}
\usepackage{color}
\nonstopmode \numberwithin{equation}{section}
\newtheorem{thm}{Theorem}
\newtheorem{lem}{Lemma}
\newtheorem{cor}{Corollary}

\newtheorem{cl}{Claim}
\newtheorem{ca}{Case}
\newtheorem{sca}{Subcase}
\newtheorem{scl}{Subclaim}
\newtheorem{conj}{Conjecture}

\theoremstyle{definition}
\newtheorem{defn}{Definition}

\newtheorem{op}[equation]{Open Problem}
\newtheorem{ques}[equation]{Question}
\newtheorem{rem}{Remark}
\newtheorem{exam}[equation]{Example}

\newcounter {own}
\def\theown {\thesection       .\arabic{own}}

\newenvironment{pf}[1][]{%
 \vskip 3mm
 \noindent
 \ifthenelse{\equal{#1}{}}%
  {{\slshape Proof. }}%
  {{\slshape #1.} }%
 }%
{\qed\bigskip}

\newcounter{alphabet}
\newcounter{tmp}
\newenvironment{Thm}[1][]{\refstepcounter{alphabet}%
\bigskip%
\noindent%
{\bf Theorem \Alph{alphabet}}%
\ifthenelse{\equal{#1}{}}{}{ (#1)}%
{\bf .} \itshape}{\vskip 8pt}

\makeatletter
\newcommand{\Ref}[1]{\@ifundefined{r@#1}{}{\setcounter{tmp}{\ref{#1}}\Alph{tmp}}}
\makeatother

\newenvironment{Lem}[1][]{\refstepcounter{alphabet}%
\bigskip%
\noindent%
{\bf Lemma \Alph{alphabet}}%
{\bf .} \itshape}{\vskip 8pt}

\newcommand{\IC}{{\mathbb C}}
\newcommand{\ID}{{\mathbb D}}

\newcommand{\diam}{{\operatorname{diam}}}

\newcommand{\dist}{{\operatorname{dist}}}

\def\be{\begin{equation}}
\def\ee{\end{equation}}

\newcommand{\bee}{\begin{enumerate}}
\newcommand{\eee}{\end{enumerate}}

\newcommand{\blem}{\begin{lem}}
\newcommand{\elem}{\end{lem}}
\newcommand{\bthm}{\begin{thm}}
\newcommand{\ethm}{\end{thm}}
\newcommand{\bcor}{\begin{cor}}
\newcommand{\ecor}{\end{cor}}
\newcommand{\beg}{\begin{exam}}
\newcommand{\eeg}{\end{exam}}
\newcommand{\begs}{\begin{examples}}
\newcommand{\eegs}{\end{examples}}
\newcommand{\bdefe}{\begin{defn}}
\newcommand{\edefe}{\end{defn}}
\newcommand{\bprob}{\begin{prob}}
\newcommand{\eprob}{\end{prob}}
\newcommand{\bques}{\begin{ques}}
\newcommand{\eques}{\end{ques}}
\newcommand{\bei}{\begin{itemize}}
\newcommand{\eei}{\end{itemize}}
\newcommand{\bcon}{\begin{conj}}
\newcommand{\econ}{\end{conj}}
\newcommand{\bop}{\begin{op}}
\newcommand{\eop}{\end{op}}

\newcommand{\bca}{\begin{ca}}
\newcommand{\eca}{\end{ca}}
\newcommand{\bsca}{\begin{sca}}
\newcommand{\esca}{\end{sca}}

\newcommand{\bcl}{\begin{cl}}
\newcommand{\ecl}{\end{cl}}

\newcommand{\bscl}{\begin{scl}}
\newcommand{\escl}{\end{scl}}

\newcommand{\bcons}{\begin{conjs}}
\newcommand{\econs}{\end{conjs}}
\newcommand{\bprop}{\begin{propo}}
\newcommand{\eprop}{\end{propo}}
\newcommand{\br}{\begin{rem}}
\newcommand{\er}{\end{rem}}
\newcommand{\brs}{\begin{rems}}
\newcommand{\ers}{\end{rems}}
\newcommand{\bo}{\begin{obser}}
\newcommand{\eo}{\end{obser}}
\newcommand{\bos}{\begin{obsers}}
\newcommand{\eos}{\end{obsers}}
\newcommand{\bpf}{\begin{pf}}
\newcommand{\epf}{\end{pf}}
\newcommand{\ba}{\begin{array}}
\newcommand{\ea}{\end{array}}
\newcommand{\beq}{\begin{eqnarray}}
\newcommand{\beqq}{\begin{eqnarray*}}
\newcommand{\eeq}{\end{eqnarray}}
\newcommand{\eeqq}{\end{eqnarray*}}

\newcommand{\ds}{\displaystyle}

\newcounter{minutes}
\divide\time by 60
\newcounter{hours}
\multiply\time by 60 \addtocounter{minutes}{-\time}

\begin{document}

\bibliographystyle{amsplain}
\title []
{John disk and $K$-quasiconformal  harmonic mappings}

\def\thefootnote{}
\footnotetext{ \texttt{\tiny File:~\jobname .tex,
          printed: \number\day-\number\month-\number\year,
          \thehours.\ifnum\theminutes<10{0}\fi\theminutes}
} \makeatletter\def\thefootnote{\@arabic\c@footnote}\makeatother

\author{Shaolin Chen}
\address{Sh. Chen, College of Mathematics and Statistics, Hengyang Normal University, Hengyang, Hunan 421008,
People's Republic of China.} \email{mathechen@126.com}

\author{Saminathan Ponnusamy 
}
\address{S. Ponnusamy,
Indian Statistical Institute (ISI), Chennai Centre, SETS (Society
for Electronic Transactions and Security), MGR Knowledge City, CIT
Campus, Taramani, Chennai 600 113, India. }
\email{samy@isichennai.res.in, samy@iitm.ac.in}




\begin{abstract}
The main aim of this article is to establish certain relationships
between $K$-quasiconformal harmonic mappings and John disks. The
results of this article are the generalizations of the corresponding
results of Ch.~Pommerenke \cite{Po}.
\end{abstract}


\maketitle \pagestyle{myheadings}
\markboth{ Sh. Chen and S. Ponnusamy}{John disk and $K$-quasiconformal  harmonic mappings}

\section{Introduction and main results }\label{csw-sec1}

For $a\in\mathbb{C}$ and  $r>0$, we let $\ID(a,r)=\{z:\, |z-a|<r\}$
so that $\mathbb{D}_r:=\mathbb{D}(0,r)$ and thus,
$\mathbb{D}:=\ID_1$ denotes the open unit disk in the complex plane
$\mathbb{C}$. This paper provides a necessary and sufficient
condition for the image $\Omega =f(\ID)$ of univalent harmonic
mappings $f$ defined on $\ID$ to be a {\it John disk } (see Theorems
\ref{cor-1} and \ref{thm-2.0} ). Some differential properties of
$K$-quasiconformal harmonic mappings will also be characterized by
using {\it Pommerenke interior domains} and John disks (see Theorem
\ref{thm-1.0} and Corollary \ref{cor-1.0}). In addition, we present
a sufficient condition, in terms of harmonic analog of the {\it
pre-Schwarzian}  of  $K$-quasiconformal harmonic mappings $f$ on
$\ID$, for $\Omega =f(\ID)$ to be a John disk (see Theorem
\ref{thm-2}). 
Similar results for
analytic functions are proved earlier by Ahlfors and Weill
\cite{AW},  Becker and Pommerenke \cite{BP}, and Pommerenke
\cite{Po}. In order to state and prove our main results and related
investigations, we need to recall some basic definitions, remarks
and some results.

For a real $2\times2$ matrix $A$, we use the matrix norm $\|A\|=\sup\{|Az|:\,|z|=1\}$ and the matrix function
$l(A)=\inf\{|Az|:\,|z|=1\}$. For $z=x+iy\in\mathbb{C}$, the formal derivative of the complex-valued functions
$f=u+iv$ is given by
$$D_{f}=\left(\begin{array}{cccc}
\ds u_{x}\;~~ u_{y}\\[2mm]
\ds v_{x}\;~~ v_{y}
\end{array}\right),
$$
so that
$$\|D_{f}\|=|f_{z}|+|f_{\overline{z}}| ~\mbox{ and }~ l(D_{f})=\big| |f_{z}|-|f_{\overline{z}}|\big |,
$$
where $f_{z}=(1/2)\big( f_x-if_y\big)$ and $f_{\overline{z}}=(1/2)\big(f_x+if_y\big)$.

Let $\Omega$ be a domain in $\mathbb{C}$, with non-empty boundary. A sense-preserving
homeomorphism $f$ from a domain $\Omega$ onto $\Omega'$, contained
in the Sobolev class $W_{loc}^{1,2}(\Omega)$, is said to be a {\it
$K$-quasiconformal mapping} if, for $z\in\Omega$,
$$\|D_{f}(z)\|^{2}\leq K\det D_{f}(z),~\mbox{i.e.,}~\|D_{f}(z)\|\leq Kl\big(D_{f}(z)\big),
$$
where $K\geq1$ and $\det D_{f}$ is the determinant of $D_{f}$ (cf.
\cite{K,LV,Va,V}).

A complex-valued function $f$ defined in a  simply connected
subdomain $G$ of $\IC$  is called a {\it harmonic mapping} in $G$ if
and only if both the real and the imaginary parts of $f$ are real
harmonic in $G$. It is indeed a simple fact that every harmonic
mapping $f$ in $G$ admits a decomposition $f=h+\overline{g}$, where
$h$ and $g$ are analytic in $G$. If we choose the additive constant
such that $g(0)=0$, then the decomposition is unique. Since the
Jacobian $\det D_{f}$ of $f$ is given by
$$\det D_{f} :=|f_{z}|^{2}-|f_{\overline{z}}|^{2} =|h'|^2-|g'|^2,
$$
$f$ is locally univalent and sense-preserving in  $G$ if and only if
$|g'(z)|<|h'(z)|$ in $G$; or equivalently if $h'(z)\neq0$ and the
dilatation $\omega =g'/h'$ has the property that $|\omega(z)|<1$ in
$G$ (see \cite{Lewy} and also \cite{Du}).

In the recent years, the family  ${\mathcal S}_{H}$ of all sense-preserving planar harmonic univalent mappings
$f=h+\overline{g}$ in $\mathbb{D}$, with the normalization $h(0)=g(0)=0$ and
$h'(0)=1$, attracted the attention of many function theorists. This class together with a few other geometric subclasses,
originally investigated in details by \cite{Clunie-Small-84}, became instrumental in
the study of univalent harmonic mappings. See the monograph \cite{Du} and the recent survey \cite{SaRa2013} for the theory of these functions.

If the co-analytic part $g$ is identically zero in the
decomposition of $f$, then the class ${\mathcal S}_{H}$
reduces to the classical family $\mathcal S$ of all normalized
analytic univalent functions $h(z)=z+\sum_{n=2}^{\infty}a_{n}z^{n}$
in $\ID$. If ${\mathcal S}_H^{0}=\{f=h+\overline{g} \in {\mathcal S}_H: \,g'(0)=0 \} $, then
the family ${\mathcal S}_H^{0}$ is both normal and compact (see \cite{Clunie-Small-84,Du,SaRa2013}).

Let $d_{\Omega}(z)$ be the Euclidean distance from $z$ to the boundary $\partial \Omega$ of
$\Omega$. In particular, we always use $d(z)$ to denote the Euclidean distance from $z$ to the boundary
$\partial \ID$ of $\mathbb{D}.$

\begin{defn}\label{CP-de-1}
A bounded simply connected plane domain $G$ is called a {\it
$c$-John disk} for $c\geq1$ with {\it John center} $w_{0}\in G$ if
for each $w_{1}\in G$ there is a rectifiable arc $\gamma$, called a
{\it John curve}, in $G$ with end points $w_{1}$ and $w_{0}$ such
that
$$\sigma_{\ell}(w)\leq cd_{G}(w)
$$
for all $w$ on $\gamma$, where $\gamma[w_{1},w]$ is the subarc of $\gamma$
between $w_{1}$ and $w$, and $\sigma_{\ell}(w)$ is the Euclidean
length of $\gamma[w_{1},w]$ (see \cite{John,KH,NV,Po1}).
\end{defn}

\begin{rem}
If $f$ is a complex-valued and univalent mapping in $\mathbb{D}$, $G=f(\mathbb{D})$ and, for $z\in\mathbb{D}$, $\gamma=f([0,z])$ in Definition
\ref{CP-de-1}, then we call $c$-John disk as a {\it radial} $c$-John
disk, where $w_{0}=f(0)$ and $w=f(z).$  In particular, if $f$ is a
conformal mapping, then we call $c$-John disk as a {\it hyperbolic}
$c$-John disk.  It is well known that any point $w_{0}\in G$ can be chosen as John center by modifying the constant
$c$ if necessary. Moveover,  if we don't emphasize the constant $c$, we regard the $c$-John
disk as the John disk (cf. \cite{John,KH,NV}).
\end{rem}

In \cite{Po} (see also \cite[p.~97]{Po1}), Pommerenke proved that if $f$ maps $\mathbb{D}$ conformally onto a bounded domain $G$, then
$G$ is a hyperbolic John disk if and only if there exist constants $M>0$ and $\delta\in(0,1)$
such that for each $\zeta\in\partial\mathbb{D}$, and for $0\leq r_{1}\leq r_{2}<1$, we have
$$|f'(r_{2}\zeta)|\leq M|f'(r_{1}\zeta)|\left(\frac{1-r_{2}}{1-r_{1}}\right)^{\delta-1}.
$$
%
%
%
Later, in \cite[Theorem 2.3]{KH}, Kari Hag and Per Hag  gave an
alternate proof of this result. In this paper, our first aim is to
extend this result to planar harmonic mappings.

\begin{thm}\label{cor-1}
For $K\geq1$, let $f\in{\mathcal S}_H^{0}$ be a $K$-quasiconformal
harmonic mapping from $\mathbb{D}$ onto a bounded domain $\Omega$.
Then $\Omega$ is a radial John disk if and only if there are
constants $M(K)>0$ and $\delta\in(0,1)$ such that for each
$\zeta\in\partial\mathbb{D}$ and for $0\leq r\leq\rho<1,$
\be\label{eq-h}\|D_{f}(\rho\zeta)\|\leq M(K)
\|D_{f}(r\zeta)\|\left(\frac{1-\rho}{1-r}\right)^{\delta-1}.
\ee
\end{thm}


The following result is another characterization of radial John
disk, which is also a generalization of \cite[Theorem 1]{Po}.

\begin{thm}\label{thm-2.0}
For $K\geq1$, let $f\in{\mathcal S}_H^{0}$ be a $K$-quasiconformal
mapping and $\Omega=f(\mathbb{D})$ is a bounded domain. Then the
following conditions are equivalent:
\bee
\item[{\rm (a)}] $\Omega$ is a radial John disk;
\item[{\rm (b)}] There is a positive constant $M_{1}$ such that, for all $z\in\mathbb{D}$,
$$
\diam f(B(z))\leq
M_{1}d_{\Omega}(f(z)),
$$
 where  $B(z)=\{
\zeta:~|z|\leq|\zeta|<1,~|\arg z-\arg \zeta|\leq\pi(1-|z|)\};$
\item[{\rm (c)}] There is a positive constant $\delta\in(0,1)$ such that, for all $z\in\mathbb{D}$ and  $\zeta\in B(z)$,
\be\label{eq-2.1}\|D_{f}(\zeta)\|\leq
M_{2}\|D_{f}(z)\|\left(\frac{1-|\zeta|}{1-|z|}\right)^{\delta-1},
\ee
where $M_{2}$ is a positive constant.
\eee
\end{thm}

By using some distortion conditions in  Theorem \ref{thm-2.0}, we
get a characterization of coefficients of $K$-quasiconformal
harmonic mappings.

\begin{thm}\label{thm-3}
For $K\geq1$, let $f=h+\overline{g}\in{\mathcal S}_H^{0}$ be a
$K$-quasiconformal harmonic mapping, where
$$h(z)=z+\sum_{n=2}^{\infty}a_{n}z^{n}~\mbox{and}~ g(z)=\sum_{n=2}^{\infty}b_{n}z^{n}.
$$
If $f$ satisfies the condition ${\rm (b)}$ or ${\rm (c)}$ in Theorem {\rm \ref{thm-2.0}}, then there is some
$\beta_{0}>0$ such that
$$\sum_{n=2}^{\infty}n^{1+\beta_{0}}(|a_{n}|^{2}+|b_{n}|^{2})<\infty.
$$
\end{thm}
%


Using Theorems \ref{thm-2.0} and \ref{thm-3}, it can be easily seen that the conclusion of Theorem~\ref{thm-3} continues to hold
if the assumption that ``$f$ satisfies the condition ${\rm (b)}$ or ${\rm (c)}$ in Theorem {\rm \ref{thm-2.0}}"
is replaced by ``$\Omega=f(\mathbb{D})$ is a radial John disk".

  For $p\in(0,\infty]$, the {\it generalized Hardy space
$H^{p}_{g}(\mathbb{D})$} consists of all those functions $f:\
\mathbb{D}\rightarrow\mathbb{C}$ such that $f$ is measurable,
$M_{p}(r,f)$ exists for all $r\in(0,1)$ and  $ \|f\|_{p}<\infty$,
where
$$\|f\|_{p}=
\begin{cases}
\displaystyle\sup_{0<r<1}M_{p}(r,f)
& \mbox{if } p\in(0,\infty)\\
\displaystyle\sup_{z\in\mathbb{D}}|f(z)| &\mbox{if } p=\infty
\end{cases},
~\mbox{ and }~
M_{p}^{p}(r,f)=\frac{1}{2\pi}\int_{0}^{2\pi}|f(re^{i\theta})|^{p}\,d\theta.
$$

Let $f\in{\mathcal S}_{H}$ be a $K$-quasiconformal harmonic mapping
from $\mathbb{D}$ onto a domain $G$. For $0<r<1$ and $w_{1},
w_{2}\in f\left(\partial\mathbb{D}_{r}\right)$, let $\gamma_{r}$ be
the smaller subarc of $f\left(\partial\mathbb{D}_{r}\right)$ between
$w_{1}$ and $w_{2}$, and let
$$d_{G_{r}}(w_{1},w_{2})=\inf_{\Gamma}\diam \Gamma,
$$
where $\Gamma$ runs through all arcs from $w_{1}$ to $w_{2}$ that
lie in $G_{r}=f(\mathbb{D}_{r})$ except for their endpoints. If
\be\label{eq-y}
\sup_{0<r<1}\left\{\sup_{w_{1},w_{2}\in\gamma_{r}}\frac{\ell\big(\gamma_{r}[w_{1},w_{2}]\big)}{d_{G_{r}}(w_{1},w_{2})}\right\}<\infty,
\ee
then we call $G$ as a {\it Pommerenke interior} domain (cf.
\cite{Po}). In particular, if $G$ is bounded, then we call $G$ as a
bounded Pommerenke interior domain. Our next theorem  is an
analogous result of  \cite[Theorem 3]{Po}.

\begin{thm}\label{thm-1.0}
Let $f\in{\mathcal S}_{H}$ be a $K$-quasiconformal harmonic mapping
from $\mathbb{D}$ onto a bounded Pommerenke interior domain $G$. If
 there are constants $M$ and $\delta\in(0,1)$ such that for each
$\zeta\in\partial\mathbb{D}$ and for $0\leq r\leq\rho<1,$
\be\label{eq-1-4}
\|D_{f}(\rho\zeta)\|\leq M \|D_{f}(r\zeta)\|\left(\frac{1-\rho}{1-r}\right)^{\delta-1},
\ee
then $\|D_{f}\|\in H^{1}_{g}(\mathbb{D}).$ 
\end{thm}

The following result easily follows from Theorems \ref{cor-1} and
\ref{thm-1.0}.

\begin{cor}\label{cor-1.0}
Let $f\in{\mathcal S}_{H}^{0}$ be a $K$-quasiconformal harmonic
mapping from $\mathbb{D}$ onto a bounded Pommerenke interior domain
$G$. If $G$ is a radial John disk, then $\|D_{f}\|\in
H^{1}_{g}(\mathbb{D}).$
\end{cor}

In terms of the canonical representation of a sense-preserving harmonic mappings $f=h+\overline{g}$ in $\mathbb{D}$ with $\omega
=g'/h'$, as in the works of Hern\'andez and Mart\'in \cite{HM}, the Pre-Schwarzian derivative $P_{f}$ of $f$
and the Schwarzian derivative $S_{f}$ of $f$ are defined by
$$ P_{f}=T_h -\frac{\omega'\overline{\omega}}{1-|\omega|^{2}},
~\mbox{ and } S_{f}=Sh+\frac{\overline{\omega}}{1-|\omega|^{2}}\left(T_h\omega'-\omega''\right)
-\frac{3}{2}\left(\frac{\omega'\overline{\omega}}{1-|\omega|^{2}}\right)^{2},
$$
respectively. Here
$$T_h= \frac{h''}{h'} ~\mbox{ and }~ Sh=T_h'-\frac{1}{2}T_h^{2}
$$
are referred to as the Pre-Schwarzian and Schwarzian (derivatives) of a locally univalent analytic function $f$ in $\ID$,
respectively.  For the original definition of the Schwarzian derivative of harmonic mappings, see \cite{CDO}.

Ahlfors and Weill \cite{AW}, Becker and Pommerenke \cite{BP}
characterized the quasidisk by using the Pre-Schwarzian of analytic
functions. On the basis of the works of Chuaqui, et al. \cite{COC},
Kari Hag and Per Hag \cite{KH} discussed the relationships between
the John disk and the Pre-Schwarzian of analytic functions. It is
natural to ask whether a similar relationship is attainable
(see \cite[Theorem 4]{COC} and \cite[Theorem 3.7]{KH}) with the help of Pre-Schwarzian of harmonic mappings. 
This is the content of our next result.

\begin{thm}\label{thm-2}
Suppose that $f\in{\mathcal S}_H^{0}$ is a $K$-quasiconformal harmonic mapping of $\mathbb{D}$ onto a bounded domain
$f(\mathbb{D})$ for some $K\geq1$ and such that
$$\lim_{|z|\rightarrow1^{-}}\sup\left\{(1-|z|^{2}){\rm Re}\big(zP_{f}(z)\big)\right\}<1.
$$
If $\ell\big(f([0,z])\big)<\infty$  for all $z\in\mathbb{D}$, then $f(\mathbb{D})$ is a radial John disk.
\end{thm}

The proofs of  Theorems \ref{cor-1}$-$\ref{thm-2}  will be given in
Section \ref{csw-sec2}.

\section{The proofs of the main results }\label{csw-sec2}

We begin the section by recalling the following results which play
an important role in the proofs of Theorems
\ref{cor-1}$-$\ref{thm-2}.

\begin{Thm} {\rm (\cite[Proposition 3.1]{MM} and \cite[Theorem 3.2]{MM})  }\label{ThmA}
Let $f$ be a $K$-quasiconformal harmonic mapping from $\mathbb{D}$
onto itself. Then for all $z\in\mathbb{D}$, we have
$$\frac{1+K}{2K}\left (\frac{1-|f(z)|^{2}}{1-|z|^{2}}\right )\leq|f_{z}(z)|\leq\frac{K+1}{2}\left (\frac{1-|f(z)|^{2}}{1-|z|^{2}}\right ).
$$
\end{Thm}

\begin{Thm} $($\cite[Theorem 3]{CPRW}$)$\label{ThmB}
Let $f\in{\mathcal S}_H^{0}$. Then there is a positive constant $c_{1}<+\infty$ such that for $\xi\in\partial\mathbb{D}$ and
$0\leq r_{3}\leq r_{4}<1$,
$$\|D_{f}(r_{4}\xi)\|\geq\frac{1}{2^{1+c_{1}}}\|D_{f}(r_{3}\xi)\|\left(\frac{1-r_{4}}{1-r_{3}}\right)^{c_{1}-1}.
$$
\end{Thm}

\subsection*{Proof of Theorem \ref{cor-1}}

We first prove the sufficiency. Applying \cite[Proposition 13]{Mi}, we obtain that
\beq\label{eq-2.0}
\|D_{f}(z)\|\leq\frac{16Kd_{\Omega}(f(z))}{1-|z|^{2}}.
\eeq
Also, by (\ref{eq-h}) and (\ref{eq-2.0}), for $w=f(r\zeta)$ and
$w_{1}=f(\rho\zeta)$, we have
\beq\label{eq-7} \nonumber
\sigma_{\ell}(w)&=&\int_{r}^{\rho}|df(t\zeta)|\leq
\int_{r}^{\rho}\|D_{f}(t\zeta)\|dt\\
\nonumber
&\leq&M(K)\|D_{f}(r\zeta)\|\int_{r}^{1}\left(\frac{1-t}{1-r}\right)^{\delta-1}dt,
~\mbox{ by \eqref{eq-h},}\\ \nonumber
&=&\frac{M(K)}{\delta}\|D_{f}(r\zeta)\|(1-r)\\ \nonumber
&\leq&\frac{M(K)}{\delta}\|D_{f}(r\zeta)\|(1-r^{2})\\
\nonumber &\leq&\frac{16KM(K)}{\delta}d_{\Omega}(w), ~\mbox{ by
\eqref{eq-2.0}},
\eeq
which implies that $\Omega$ is a radial $(16KM(K)/\delta)$-John disk with John center $w_{0}=f(0)$ and with
$\gamma=f([0,\rho\zeta])$ as the John curves, where $r\in[0,1)$, $\rho\in[r,1)$ and $\zeta\in\partial\mathbb{D}$.

Now we  prove the necessity. For $z\in\mathbb{D}$, let
$$\Delta=f^{-1}\Big(\mathbb{D}\big(f(z), d_{\Omega}(f(z))\big)\Big)
$$
and $\phi$ be a conformal mapping of $\mathbb{D}$ onto $\Delta$ with $\phi(0)=z$. Since $\phi(\mathbb{D})\subset\mathbb{D}$,
we know that, for $w\in\mathbb{D}$,
\be\label{eq-Schw}
|\phi'(w)|\leq\frac{1-|\phi(w)|^{2}}{1-|w|^{2}}.
\ee
Then
$$F(w)=\frac{1}{d_{\Omega}(f(z))}\big(f(\phi(w))-f(z)\big)
$$
is a $K$-quasiconformal harmonic mapping of $\mathbb{D}$ onto itself
with $F(0)=0.$ It is not difficult to know that
$$\|D_{F}(w)\|=\frac{|\phi'(w)|\|D_{f}(\phi(w))\|}{d_{\Omega}(f(z))},
$$
which, together with (\ref{eq-Schw}) and Theorem \Ref{ThmA}, give that
\beq\label{eq-2}
\nonumber
\|D_{f}(z)\| & =& \|D_{f}(\phi(0))\|=
\frac{d_{\Omega}(f(z))\|D_{F}(0)\|}{|\phi'(0)|}\\
\nonumber
&\geq&\frac{d_{\Omega}(f(z))\|D_{F}(0)\|}{1-|z|^{2}}\\
&\geq&\frac{1+K}{2K}\frac{d_{\Omega}(f(z))}{1-|z|^{2}}.
\eeq
Since $\Omega$ is a radial John disk, we can choose $w_{0}=f(0)$ as the
John center and $\gamma=f([0,\rho\zeta])$ as the John curve;
$\Omega$ can be assumed to be a radial $c$-John disk with respect to
this choice, where $c\geq 1$. Hence for $w=f(r\zeta)$ and
$w_{1}=f(\rho\zeta)$, we have
\be\label{eq-3y1}
\sigma_{\ell}(w)\leq cd_{\Omega}(w)~\mbox{for all}~\rho\in[r,1).
\ee
The boundedness of $\Omega$ implies that $d_{\Omega}(w)$ is finite for all
$w\in\Omega$. Hence the limit
\be\label{eq-3y2}
\lim_{\rho\rightarrow1-}\int_{r}^{\rho}|df(t\zeta)|
\ee
does exist and is finite. By (\ref{eq-3y1}) and (\ref{eq-3y2}), we get
\be\label{eq-3}
\frac{1}{K}\int_{r}^{1}\|D_{f}(t\zeta)\|dt\leq\int_{r}^{1}l(D_{f}(t\zeta))dt\leq\int_{r}^{1}|df(t\zeta)|
\leq cd_{\Omega}(w),
\ee
where $\zeta\in\partial\mathbb{D}.$ By (\ref{eq-2}) and (\ref{eq-3}), we have
\be\label{eq-4}
\int_{r}^{1}\|D_{f}(t\zeta)\|dt\leq\frac{2cK^{2}}{1+K}(1-r^{2})\|D_{f}(r\zeta)\|
\leq M_{0}(1-r)\|D_{f}(r\zeta)\|,
\ee
where $M_{0}=\frac{4cK^{2}}{1+K}\geq2c.$

Next, we let
$$\varphi(r)=(1-r)^{-\frac{1}{M_{0}}}\int_{r}^{1}\|D_{f}(t\zeta)\|dt.
$$
By (\ref{eq-4}), we have
$$\varphi'(r)=(1-r)^{-\frac{1}{M_{0}}}\left[\frac{1}{M_{0}(1-r)}\int_{r}^{1}\|D_{f}(t\zeta)\|dt-\|D_{f}(r\zeta)\|\right]\leq0,
$$
which implies that $\varphi (r)$ is decreasing on the unit interval $(0,1).$

By Theorem \Ref{ThmB}, for $\rho\leq t\leq\frac{1+\rho}{2},$ there
is a positive constant $c_{1}$ such that
$$\|D_{f}(\rho\zeta)\|\leq 4^{c_{1}}\|D_{f}(t\zeta)\|,
$$
which gives  
\beq\label{eq-5} \nonumber
\int_{\rho}^{1}\|D_{f}(t\zeta)\|dt&\geq&\int_{\rho}^{\frac{1+\rho}{2}}\|D_{f}(t\zeta)\|dt\\
\nonumber &\geq&
4^{-c_{1}}\|D_{f}(\rho\zeta)\|\int_{\rho}^{\frac{1+\rho}{2}}dt\\
&=&2^{-2c_{1}-1}\|D_{f}(\rho\zeta)\|(1-\rho).
\eeq
For $0\leq r\leq\rho<1$, by (\ref{eq-4}) and (\ref{eq-5}), we have
\beq\label{eq-6} \nonumber
(1-\rho)^{1-\frac{1}{M_{0}}}\|D_{f}(\rho\zeta)\|&\leq&2^{1+2c_{1}}\varphi(\rho)\leq2^{1+2c_{1}}\varphi(r)\\
\nonumber
&\leq&2^{1+2c_{1}}M_{0}(1-r)^{1-\frac{1}{M_{0}}}\|D_{f}(r\zeta)\|,
\eeq
which yields
\begin{eqnarray*}
\|D_{f}(\rho\zeta)\|&\leq&2^{1+2c_{1}}M_{0}\|D_{f}(r\zeta)\|\left(\frac{1-r}{1-\rho}\right)^{1-\frac{1}{M_{0}}}\\
&=&2^{1+2c_{1}}M_{0}\|D_{f}(r\zeta)\|\left(\frac{1-\rho}{1-r}\right)^{\frac{1}{M_{0}}-1}.
\end{eqnarray*}
The proof of the theorem is complete.  \hfill $\Box$

\bigskip

For $z_{1}, z_{2}\in\mathbb{D}$, the {\it hyperbolic metric} (or
{\it Poincar\'e metric}) is defined by

$$\lambda_{\mathbb{D}}(z_{1}, z_{2})=\min_{\gamma}\int_{\gamma}\frac{|dz|}{1-|z|^{2}},
$$
where the minimum is taken over all curves $\gamma$ in $\mathbb{D}$ from $z_{1}$ and $z_{2}$. It is well-known that, for
$z_{1}, z_{2}\in\mathbb{D}$,
$$\lambda_{\mathbb{D}}(z_{1}, z_{2})=\frac{1}{2}\log\frac{1+|z_{1}-z_{2}|/|1-\overline{z}_{1}z_{2}|}{1-|z_{1}-z_{2}|/|1-\overline{z}_{1}z_{2}|},
$$
which is equivalent to
$$\left|\frac{z_{1}-z_{2}}{1-\overline{z}_{1}z_{2}}\right|=\frac{e^{2\lambda_{\mathbb{D}}(z_{1},
z_{2})}-1}{e^{2\lambda_{\mathbb{D}}(z_{1},
z_{2})}+1}=\tanh\lambda_{\mathbb{D}}(z_{1}, z_{2}).
$$

In \cite{Small}, Sheil-Small proved the following result.

\begin{Lem}\label{Lem-B}
Let $f=h+\overline{g}\in{\mathcal S}_{H}$ and
$\alpha=\sup_{f\in{\mathcal S}_H} \frac{|h''(0)|}{2}$, where $h$
and $g$ are analytic in $\mathbb{D}$. Then
$$\frac{(1-|z|)^{\alpha-1}}{(1+|z|)^{\alpha+1}}\leq|h'(z)|\leq\frac{(1+|z|)^{\alpha-1}}{(1-|z|)^{\alpha+1}}.
$$
\end{Lem}

We remark that $\alpha=\sup_{f\in{\mathcal
S}_{H}}\frac{|h''(0)|}{2}$ is finite, but the sharp upper bound of
$\alpha$ is still unknown (see \cite{Du,Small}).

\begin{lem}\label{lem-1}
Let $f=h+\overline{g}\in{\mathcal S}_{H}$, where $h$ and $g$ are
analytic in $\mathbb{D}$. Then, for $z_{1}, z_{2}\in\mathbb{D}$,
$$\frac{1}{2}\|D_{f}(z_{1})\|e^{-(1+\alpha)\lambda_{\mathbb{D}}(z_{1},
z_{2})}\leq\|D_{f}(z_{2})\|\leq2\|D_{f}(z_{1})\|e^{(1+\alpha)\lambda_{\mathbb{D}}(z_{1}, z_{2})},
$$
where $\alpha$ is defined in Lemma \Ref{Lem-B}.
\end{lem}

\bpf Let $f=h+\overline{g}\in{\mathcal S}_{H}$ and
$z=\frac{z_{2}-z_{1}}{1-\overline{z}_{1}z_{2}}$, where $h$, $g$ are
analytic in $\mathbb{D}$ and $z_{1}, z_{2}\in\mathbb{D}$. Then
$$F(z)=\frac{f(z_{2})-f(z_{1})}{(1-|z_{1}|^{2})h'(z_{1})}\in{\mathcal S}_{H},
$$
where $z_{2}=\frac{z+z_{1}}{1+\overline{z}_{1}z}.$ By Lemma \Ref{Lem-B}, we get
$$\frac{(1-|z|)^{\alpha-1}}{(1+|z|)^{\alpha+1}}\leq|F_{z}(z)|=\frac{|h'(z_{2})|}{|h'(z_{1})||1+\overline{z}_{1}z|^{2}}
\leq\frac{(1+|z|)^{\alpha-1}}{(1-|z|)^{\alpha+1}},
$$
which gives
\be\label{eq-1a1}
\frac{(1-|z|)^{\alpha+1}}{(1+|z|)^{\alpha+1}}|h'(z_{1})|\leq|h'(z_{2})|
\leq\frac{(1+|z|)^{\alpha+1}}{(1-|z|)^{\alpha+1}}|h'(z_{1})|.
\ee
By (\ref{eq-1a1}), we obtain
$$\frac{1}{2}\frac{(1-|z|)^{\alpha+1}}{(1+|z|)^{\alpha+1}}\|D_{f}(z_{1})\|\leq\|D_{f}(z_{2})\|\leq
2\frac{(1+|z|)^{\alpha+1}}{(1-|z|)^{\alpha+1}}\|D_{f}(z_{1}),
$$
which implies that
$$\frac{1}{2}\|D_{f}(z_{1})\|e^{-(1+\alpha)\lambda_{\mathbb{D}}(z_{1},
z_{2})}\leq\|D_{f}(z_{2})\|\leq2\|D_{f}(z_{1})\|e^{(1+\alpha)\lambda_{\mathbb{D}}(z_{1}, z_{2})}.
$$
The proof of this lemma is complete. \epf

\begin{lem}\label{lem-2}
Let $a_{1}, a_{2}$ and $a_{3}$ be positive constants and let
$0<|z_{0}|=1-\delta$, where $\delta\in(0,1)$. If $f\in{mathcal
S}_{H}$, $0\leq1-a_{2}\delta\leq|z|\leq1-a_{1}\delta$ and $|\arg z-
\arg z_{0}|\leq a_{3}\delta$, then
$$\frac{1}{M(a_{1},a_{2},a_{3})}\|D_{f}(z_{0})\|\leq\|D_{f}(z)\|\leq
M(a_{1},a_{2},a_{3})\|D_{f}(z_{0})\|,
$$
where $M(a_{1},a_{2},a_{3})=2e^{(1+\alpha)\left(a_{3}+\frac{1}{2}\log\frac{2a_{2}-a_{1}}{a_{1}}\right)}$
and $\alpha$ is defined in Lemma \Ref{Lem-B}.
\end{lem}

\begin{figure}[!ht]
\centering
\includegraphics{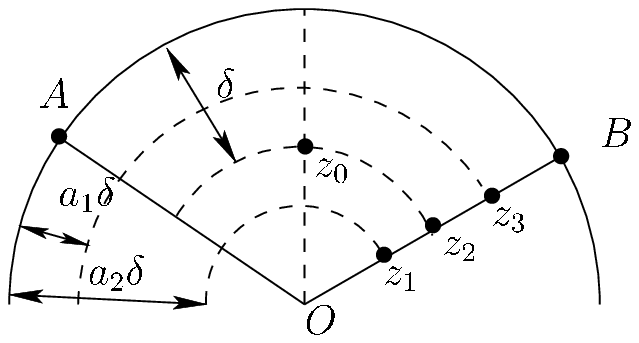}
\caption{ }
\label{fig1}
\end{figure}

\bpf Let $\angle AOB=2a_{3}\delta$ and $z_{1}, z_{2}, z_{3}$ line in
the line $OB$ with $|z_{1}|\leq|z_{2}|=|z_{0}|\leq|z_{3}|$, see Figure \ref{fig1}.
 Then the
length of the circular arc from $z_{0}$ to $z_{2}$ is less than
$a_{3}\delta$. By calculations, we have
$$\lambda_{\mathbb{D}}(z_{0}, z_{2})<\frac{a_{3}\delta}{1-(1-\delta)^{2}}=\frac{a_{3}}{2-\delta}<a_{3}
$$
and
$$\left|\frac{z_{3}-z_{1}}{1-\overline{z}_{1}z_{3}}\right|=\frac{1-a_{1}\delta-(1-a_{2}\delta)}
{1-(1-a_{1}\delta)(1-a_{2}\delta)}=\frac{a_{2}-a_{1}}{a_{2}+a_{1}(1-a_{2}\delta)}\leq\frac{a_{2}-a_{1}}{a_{2}}.
$$
Hence
\begin{eqnarray*}
\lambda_{\mathbb{D}}(z_{0}, z)&\leq&\lambda_{\mathbb{D}}(z_{0},
z_{2})+\lambda_{\mathbb{D}}(z_{2}, z_{1})\\
&\leq&\lambda_{\mathbb{D}}(z_{0}, z_{2})+\lambda_{\mathbb{D}}(z_{1},
z_{3})\\
&\leq&a_{3}+\frac{1}{2}\log\frac{2a_{2}-a_{1}}{a_{1}}.
\end{eqnarray*}
By Lemma \ref{lem-1}, we see that
$$\frac{1}{M(a_{1},a_{2},a_{3})}\|D_{f}(z_{0})\|\leq\|D_{f}(z)\|\leq M(a_{1},a_{2},a_{3})\|D_{f}(z_{0})\|,
$$
where $M(a_{1},a_{2},a_{3})$ is defined as in the statement.
\epf

\subsection*{Proof of Theorem \ref{thm-2.0}}
We first prove (c)$\Rightarrow$(b). Let
$z=re^{i\theta}\in\mathbb{D}$ and $r_{1}e^{i\theta_{1}},
r_{2}e^{i\theta_{2}}\in B(re^{i\theta})$ with $r_{1}\leq r_{2}$.
Then, by (\ref{eq-2.1}), Lemma \ref{lem-2} and \cite[Proposition 13]{Mi}, there is a positive constant $M$ such that
\begin{eqnarray*}
|f(r_{2}e^{i\theta_{2}})-f(r_{1}e^{i\theta_{1}})|&\leq&|f(r_{2}e^{i\theta_{2}})-f(re^{i\theta_{2}})|+|f(r_{1}e^{i\theta_{1}})-f(re^{i\theta_{1}})|\\
&&+|f(re^{i\theta_{2}})-f(re^{i\theta_{1}})|\\
&\leq&\int_{r}^{r_{2}}\|D_{f}(\rho
e^{i\theta_{2}})\|d\rho+\int_{r}^{r_{1}}\|D_{f}(\rho
e^{i\theta_{1}})\|d\rho\\
&& +r\int_{\gamma_{0}}\|D_{f}(r e^{it})\|dt\\
&\leq&M_{2}\int_{r}^{r_{2}}\|D_{f}(re^{i\theta})\|\left(\frac{1-\rho}{1-r}\right)^{\delta-1}d\rho
\\&&
+M_{2}\int_{r}^{r_{1}}\|D_{f}(re^{i\theta})\|\left(\frac{1-\rho}{1-r}\right)^{\delta-1}d\rho\\
&&+Mr\int_{\gamma_{0}}\|D_{f}(r e^{i\theta})\|dt~\mbox{ (by Lemma \ref{lem-2})}\\
&\leq&\frac{2M_{2}}{\delta}\|D_{f}(re^{i\theta})\|(1-r)+Mr\ell(\gamma_{0})\|D_{f}(re^{i\theta})\|\\
&\leq&\frac{2M_{2}}{\delta}\|D_{f}(re^{i\theta})\|(1-r)+M|\theta_{2}-\theta_{1}|\|D_{f}(re^{i\theta})\|\\
&\leq&\left(\frac{2M_{2}}{\delta}+2\pi
M\right)\|D_{f}(re^{i\theta})\|(1-r)\\
&\leq&16K\left(\frac{2M_{2}}{\delta}+2\pi
M\right)d_{\Omega}(f(z)),~\mbox{ by \cite[Proposition 13]{Mi}},
\end{eqnarray*}
where $\gamma_{0}$ is the smaller subarc of $\partial\mathbb{D}_{r}$
between $re^{i\theta_{1}}$ and $re^{i\theta_{2}}$. Hence there
exists a positive constant $M_{1}$ such that, for all
$z\in\mathbb{D}$,
$$\diam f(B(z))\leq M_{1}d_{\Omega}(f(z)).
$$

Next we prove (b)$\Rightarrow$(c). For
$z=re^{i\theta}\in\mathbb{D}$, let
\be\label{eq-3cp}
\phi(r)=\int_{r}^{1}(1-x)\|D_{f}(xe^{i\theta})\|^{2}dx
\ee
and
$$\Delta(r)=\{\zeta=x+iy:~r\leq x<1,~0\leq y\leq1-x\}.
$$
By Lemma \ref{lem-2}, for $\zeta=x+iy\in\Delta(r)$, there exists a positive
constant $M_{3}$ such that
$$\|D_{f}(xe^{i\theta})\|\leq M_{3}\|D_{f}(\zeta e^{i\theta})\|,
$$
which implies that
\beq\label{eq-4cp}
\phi(r)&\leq&\int_{r}^{1}\int_{0}^{1-x}\|D_{f}(xe^{i\theta})\|^{2}dydx\\
\nonumber &\leq&M_{3}^{2}\int_{r}^{1}\int_{0}^{1-x}\|D_{f}(\zeta
e^{i\theta})\|^{2}dydx\\ \nonumber
&\leq&KM_{3}^{2}\int_{r}^{1}\int_{0}^{1-x}J_{f}(\zeta
e^{i\theta})dydx\\
\nonumber&=&KM_{3}^{2}A\big(f(\Delta(re^{i\theta}))\big),
\eeq
where
$$\Delta(re^{i\theta})=\{\zeta e^{i\theta}=(x+iy)e^{i\theta}:~r\leq x<1,~0\leq y\leq1-x\}.
$$
It is not difficult to see that $\Delta(re^{i\theta})\subset
B(re^{i\theta})$, which, together with  (\ref{eq-2}) and (\ref{eq-4cp}), imply
\beq\label{eq-5cp}
\phi(r)&\leq&KM_{3}^{2}A\big(f(\Delta(re^{i\theta}))\big)\leq KM_{3}^{2}A\big(f(B(re^{i\theta}))\big)\\
\nonumber&\leq&\frac{\pi KM_{3}^{2}}{4}\left(\diam
\big(f(B(re^{i\theta})\big)\right)^{2}\\ \nonumber&\leq&\frac{\pi
KM_{3}^{2}M_{1}^{2}}{4}\left(d_{\Omega}(f(z))\right)^{2}\\ \nonumber
&\leq&\frac{\pi
K^{3}M_{1}^{2}M_{3}^{2}}{(1+K)^{2}}(1-|z|^{2})^{2}\|D_{f}(z)\|^{2}.
\eeq By (\ref{eq-3cp}), for $r\leq\rho<1$, we get \be\label{eq-6cp}
\log\frac{\phi(\rho)}{\phi(r)}=\int_{r}^{\rho}\frac{\phi'(t)}{\phi(t)}dt\leq-\alpha\int_{r}^{\rho}\frac{dt}{1-t}
=\alpha\log\frac{1-\rho}{1-r}, \ee where $\alpha=(1+K)^{2}/\left[\pi
K^{3}M_{1}^{2}M_{3}^{2}\right].$ For $\rho\leq
x\leq\frac{1+\rho}{2}$, by Theorem \Ref{ThmB}, there is a positive
constant $c_{1}^{\ast}$ such that \be\label{eq-7cp} \|D_{f}(\rho
e^{i\theta})\|\leq 4^{c_{1}^{\ast}}\|D_{f}(x e^{i\theta})\|. \ee

Applying (\ref{eq-3cp}), (\ref{eq-6cp}) and (\ref{eq-7cp}), we have
\beq\label{eq-8cp} \nonumber
\frac{1}{2^{4c_{1}^{\ast}+1}}(1-\rho)^{2}\|D_{f}(\rho
e^{i\theta})\|^{2}&=&\frac{1}{2^{4c_{1}^{\ast}}}\int_{\rho}^{1}(1-x)\|D_{f}(\rho
e^{i\theta})\|^{2}dx\\ \nonumber &\leq&\int_{\rho}^{1}(1-x)\|D_{f}(x
e^{i\theta})\|^{2}dx=\phi(\rho)\\
\nonumber&\leq&\phi(r)\left(\frac{1-\rho}{1-r}\right)^{\alpha},
\eeq
which, together with (\ref{eq-5cp}), yield that
$$\frac{1}{2^{4c_{1}^{\ast}+1}}\|D_{f}(\rho
e^{i\theta})\|^{2}(1-\rho)^{2}\leq\phi(r)\left(\frac{1-\rho}{1-r}\right)^{\alpha}
\leq\frac{1}{\alpha}(1-r)^{2}\|D_{f}(r
e^{i\theta})\|^{2}\left(\frac{1-\rho}{1-r}\right)^{\alpha}.
$$
Then we conclude that \beq\label{eq-9cp} \|D_{f}(\rho
e^{i\theta})\|\leq\sqrt{\frac{2^{1+4c_{1}^{\ast}}}{\alpha}}\|D_{f}(z)\|\left(\frac{1-\rho}{1-r}\right)^{\frac{\alpha}{2}-1}.
\eeq By (\ref{eq-9cp}) and Lemma \ref{lem-2}, for all $\zeta=\rho
e^{i\eta}\in B(z)$, there exists a positive constant $M_{4}$ such
that
$$\|D_{f}(\rho e^{i\eta})\|\leq M_{4}\|D_{f}(\rho e^{i\theta})\|\leq M_{4}
\sqrt{\frac{2^{1+4c_{1}^{\ast}}}{\alpha}}\|D_{f}(z)\|\left(\frac{1-\rho}{1-r}\right)^{\frac{\alpha}{2}-1}.
$$

Now we  prove (a)$\Rightarrow$(c). By Theorem \ref{cor-1},
there are constants $M$ and $\delta\in(0,1)$ such that for each
$\zeta\in\partial\mathbb{D}$ and for $0\leq r\leq\rho<1,$
\be\label{eq-1-4a}
\|D_{f}(\rho\zeta)\|\leq M \|D_{f}(r\zeta)\|\left(\frac{1-\rho}{1-r}\right)^{\delta-1}.
\ee
For all $\xi\in\partial\mathbb{D}$ with $|\arg\xi-\arg\zeta|\leq\pi(1-r)$, by Lemma \ref{lem-2}, there is a
positive constant $M'$ such that
\be\label{eq-1-5}
\|D_{f}(r\zeta)\|\leq M'\|D_{f}(r\xi)\|.
\ee
Hence (\ref{eq-2.1}) follows from (\ref{eq-1-4a}) and (\ref{eq-1-5}).

At last, we  prove (c)$\Rightarrow$(a). By \cite[Proposition
13]{Mi}, (\ref{eq-2.1}), for $w=f(r\zeta)$ and $w_{1}=f(\rho\zeta)$,
we have
\begin{eqnarray*}
\sigma_{\ell}(w)&=&\int_{r}^{\rho}|df(t\zeta)|\leq\int_{r}^{\rho}\|D_{f}(t\zeta)\|dt\\
&\leq&M_{2}\|D_{f}(r\zeta)\|\int_{r}^{\rho}\left(\frac{1-t}{1-r}\right)^{\delta-1}dt
=\frac{M_{2}}{\delta}\|D_{f}(r\zeta)\|(1-r)\\
&\leq&\frac{M_{2}}{\delta}\|D_{f}(r\zeta)\|(1-r^{2})\\
&\leq&\frac{16KM_{2}}{\delta}d_{\Omega}(w),
\end{eqnarray*}
which implies that $\Omega$ is a radial $(16KM_{2}/\delta)$-John
disk with John center $0$ and with $\gamma=f([0,\rho\zeta])$ as the
John curves, where $r\in[0,1)$, $\rho\in[r,1)$ and
$\zeta\in\partial\mathbb{D}.$ The proof is complete. \hfill $\Box$

\subsection*{Proof of Proposition \ref{thm-3}} Without loss of
generality, we assume that there is a positive constant
$M_{1}^{\star}$ such that, for all $z\in\mathbb{D}$,
\be\label{eq-2m} \diam f(B(z))\leq M_{1}^{\star}d_{\Omega}(f(z)),
\ee where $\Omega=f(\mathbb{D}).$ For $r\in[0,1)$, let

\beq\label{eqy}
\varphi(r)&=&\frac{1}{2\pi}\int_{0}^{2\pi}\left(|f_{z}(re^{it})|^{2}+|f_{\overline{z}}(re^{it})|^{2}\right)dt\\
\nonumber
&=&1+\sum_{n=2}^{\infty}n^{2}\left(|a_{n}|^{2}+|b_{n}|^{2}\right)r^{2n-2}.
\eeq Then, by Theorem \Ref{ThmA} and  (\ref{eq-2m}), we obtain
\beq\label{eq-1-8} \int_{r}^{1}\int_{-\pi(1-r)}^{\pi(1-r)}J_{f}(\rho
e^{i(\theta+t)})\rho d\theta d\rho&=&A\left(f(B(re^{it}))\right)\\
\nonumber &\leq&\frac{\pi}{4}\diam^{2}\left(f(B(re^{it}))\right)\\
\nonumber&\leq&M^{\ast}(1-r^{2})^{2}\|D_{f}(re^{it})\|^{2}, \eeq
where $M^{\ast}=\frac{\pi K^{2}M_{1}^{\star 2}}{(1+K)^{2}}.$

By (\ref{eq-1-8}), for $r\in[\frac{1}{2},1)$, we obtain
\begin{eqnarray*}
\frac{1}{2K}\int_{r}^{1}\int_{-\pi(1-r)}^{\pi(1-r)}\varphi(\rho)d\theta
d\rho&\leq&\frac{1}{K}\int_{r}^{1}\rho\left(\int_{0}^{2\pi}\|D_{f}(\rho
e^{i(t+\theta)})\|^{2}dt\right)d\theta d\rho\\
&\leq&\int_{r}^{1}\int_{-\pi(1-r)}^{\pi(1-r)}\rho\left(\int_{0}^{2\pi}J_{f}(\rho
e^{i(t+\theta)})\right)d\theta d\rho\\
&\leq&4M^{\ast}(1-r)^{2}\int_{0}^{2\pi}\|D_{f}(re^{it})\|^{2}dt\\
&\leq&16\pi M^{\ast}(1-r)^{2}\varphi(r),
\end{eqnarray*}
which gives that
\be\label{eq-1-9.0}
\int_{r}^{1}\varphi(\rho)d\rho\leq16KM^{\ast}(1-r)\varphi(r)=\beta(1-r)\varphi(r),
\ee
where $\beta=16KM^{\ast}$. Applying (\ref{eq-1-9.0}), for $r\in[\frac{1}{2},1)$, we get
\beq\label{eq-1-9}
&&\frac{d}{dr}\left[(1-r)^{-2\beta_{0}}\int_{r}^{1}\varphi(\rho)d\rho\right]\\
\nonumber &=&
\frac{1}{2\beta_{0}}(1-r)^{-2\beta_{0}-1}\int_{r}^{1}\varphi(\rho)d\rho-(1-r)^{-2\beta_{0}}\varphi(r) \leq0,
\eeq
where $\beta_{0}=1/(2\beta)$. By (\ref{eq-1-9}), for $r\in[\frac{1}{2},1)$, we have
\be\label{eq-1-10}
(1-r)^{1-2\beta_{0}}\varphi(r)\leq(1-r)^{-2\beta_{0}}\int_{r}^{1}\varphi(\rho)
d\rho\leq2^{-2\beta_{0}}\int_{\frac{1}{2}}^{1}\varphi(\rho)d\rho<\infty.
\ee

It follows from (\ref{eqy}) and (\ref{eq-1-10})  that there are two
positive constants $M_{1}^{'}$ and $M_{1}^{''}$ such that
\begin{eqnarray*}
1+\sum_{n=2}^{\infty}n^{1+\beta_{0}}(|a_{n}|^{2}+|b_{n}|^{2})&\leq&M_{1}^{'}\int_{\frac{1}{2}}^{1}(1-r)^{-\beta_{0}}\varphi(r)dr\\
&\leq&M_{1}^{''}\int_{\frac{1}{2}}^{1}(1-r)^{\beta_{0}-1}dr<\infty.
\end{eqnarray*}
The proof of this proposition is complete. \hfill $\Box$

\begin{lem}\label{lem-2p}
Let $f\in{\mathcal S}_{H}$ be a $K$-quasiconformal harmonic mapping
from $\mathbb{D}$ onto a bounded domain $G$. If there are constants
$M$ and $\delta\in(0,1)$ such that for each $\varsigma\in\partial\mathbb{D}$ and for
$0\leq r\leq\rho<1,$
\be\label{eq-14}
\|D_{f}(\rho\varsigma)\|\leq M
\|D_{f}(r\varsigma)\|\left(\frac{1-\rho}{1-r}\right)^{\delta-1},
\ee
then, for $a\in\mathbb{D}$, we have
$$\diam f(I(a))\leq M_{0}'d_{G}(a),
$$
where
$$I(a)=\{z\in\partial\mathbb{D}:\, |\arg z-\arg a|\leq1-|a|\}
$$
and
$$M_{0}'=32K\left(2e^{(1+\alpha)} +\frac{M2e^{(1+\alpha)}}{\delta}+\frac{M}{\delta}\right).
$$

\end{lem}
\bpf For $a\in\mathbb{D}$, let $a=\rho\zeta$ with $\rho=|a|$. For
$z\in I(a)$, by (\ref{eq-14}) and Lemma \ref{lem-2}, we have
\beq\label{eq-15}
|f(z\rho)-f(\rho\zeta)|&\leq&\int_{\gamma'}\rho\|D_{f}(\rho\xi)\| \,|d\xi|\\
\nonumber
&\leq&2e^{(1+\alpha)}\rho\int_{\gamma'}\|D_{f}(\rho\zeta)\|
\,|d\xi|, ~\mbox{ by Lemma \ref{lem-2},}\\
\nonumber&=&2e^{(1+\alpha)}\rho\ell(\gamma')\|D_{f}(\rho\zeta)\|\\
\nonumber
&=&2e^{(1+\alpha)}\rho^{2}\|D_{f}(\rho\zeta)\| \,|\arg (\rho\zeta)-\arg z|\\
\nonumber&\leq&2e^{(1+\alpha)}\rho^{2}(1-\rho)\|D_{f}(\rho\zeta)\|\\
\nonumber &\leq&2e^{(1+\alpha)}(1-\rho)\|D_{f}(\rho\zeta)\|,
\eeq
\beq\label{eq-16}
|f(z\rho)-f(z)|&\leq&\int_{\rho}^{1}\|D_{f}(tz)\|\,dt  \\
\nonumber
&\leq&M\int_{\rho}^{1}\|D_{f}(\rho z)\|\left(\frac{1-t}{1-\rho}\right)^{\delta-1}dt, ~\mbox{by (\ref{eq-14}),}\\
\nonumber &=&\frac{M}{\delta}(1-\rho)\|D_{f}(\rho z)\|\\ \nonumber
&\leq&\frac{2Me^{(1+\alpha)}}{\delta}(1-\rho)\|D_{f}(\rho \zeta)\|
\eeq
and
\beq\label{eq-17}
|f(\zeta\rho)-f(\zeta)|&\leq&\int_{\rho}^{1}\|D_{f}(t\zeta)\|dt  \\
\nonumber
&\leq&M\int_{\rho}^{1}\|D_{f}(\rho\zeta)\|\left(\frac{1-t}{1-\rho}\right)^{\delta-1}dt, ~\mbox{by (\ref{eq-14}),}\\
\nonumber &=&\frac{M}{\delta}(1-\rho)\|D_{f}(\rho\zeta)\|,
\eeq
where $\gamma'$ is the smaller subarc of $\partial\mathbb{D}_{\rho}$
between $\rho z$ and $\rho\zeta$.

Again, for $z\in I(a)$, by (\ref{eq-2.0}), (\ref{eq-15}), (\ref{eq-16}) and (\ref{eq-17}), we obtain
\beq
\nonumber |f(\zeta)-f(z)|&\leq& |f(\rho\zeta)-f(\rho z)|+|f(z)-f(\rho z)|+|f(\rho\zeta)-f(\zeta)|\\
\nonumber &\leq&M_{1}^{\ast}(1-\rho)\|D_{f}(\rho\zeta)\|\\
\nonumber &\leq&16M_{1}^{\ast}Kd_{G}(a),  ~\mbox{by (\ref{eq-2.0}),}
\eeq
which in turn implies that $\diam f(I(a))\leq 32KM_{1}^{\ast}d_{G}(a), $ where
\be\label{eqex1d}
M_{1}^{\ast}=2e^{(1+\alpha)} +\frac{M2e^{(1+\alpha)}}{\delta}+\frac{M}{\delta}.
\ee
The proof of the lemma is complete. \epf

\subsection*{Proof of Theorem \ref{thm-1.0}}  Let $\frac{1}{2}<\nu<1$ and
\be\label{eq-22}
\sup_{0<r<1}\left\{\sup_{w_{1},w_{2}\in\gamma_{r}}\frac{\ell\big(\gamma_{r}
[w_{1},w_{2}]\big)}{d_{G_{r}}(w_{1},w_{2})}\right\}=M_{\gamma},
\ee
where $\gamma_{r}$ is given by (\ref{eq-y}). Then, by (\ref{eq-22}),
Lemma \ref{lem-2p} and \cite[Theorem 3]{Co}, we have
\beq \nonumber
\frac{\nu}{K}\int_{0}^{2\pi}\|D_{f}(\nu
e^{i\theta})\|d\theta&\leq&\nu\int_{0}^{2\pi}l\big(D_{f}(\nu
e^{i\theta})\big)d\theta\\
\nonumber
&\leq&\int_{0}^{2\pi}|d f(\nu e^{i\theta})|\\
\nonumber&\leq&\sum_{k=1}^{7}\int_{I(z_{k})}|d
f(\nu e^{i\theta})|\\
\nonumber&\leq&M_{\gamma}\sum_{k=1}^{7}\diam f\left(I(z_{k})\right),
~\mbox{by (\ref{eq-22}),}
\\
\nonumber&\leq&32M_{\gamma}M_{1}^{\ast}K\sum_{k=1}^{7}d_{G}\left(f(z_{k})\right),~\mbox{by Lemma \ref{lem-2p},}\\
\nonumber
&\leq&\frac{64M_{\gamma}M_{1}^{\ast}K}{1+K}\sum_{k=1}^{7}\left\{(1-|z_{k}|^{2})\|D_{f}(z_{k})\|\right\}\\
\nonumber &\leq&\frac{1792M_{\gamma}M_{1}^{\ast}K}{(1+K)\pi},
~\mbox{by \cite[Theorem 3]{Co},}
\eeq
which implies that $\|D_{f}\|\in H^{1}_{g}(\mathbb{D})$, where $k\in\{1,2,...,7\}$,
$$z_{k}=\frac{1}{2}e^{i(k-1)},~~ I(z_{k})=\{z\in\partial\mathbb{D}:\, |\arg z-\arg
z_{k}|\leq1-|z_{k}|\} ,
$$
and $M_{1}^{\ast}$ is given by \eqref{eqex1d}. The proof of the
theorem is complete. \hfill $\Box$

\subsection*{Proof of Theorem \ref{thm-2}}
 By the assumption, we see that there is a $\nu\in(0,1)$ and
$r_{0}\in(0,1)$ such that, for $r_{0}\leq\eta<1,$
$$\frac{\nu}{1-\eta^{2}}\geq{\rm Re}\big(\zeta P_{f}(\eta\zeta)\big)=\mbox{Re}\left(\frac{\zeta h''(\eta\zeta)}{h'(\eta\zeta)}\right)-
\mbox{Re}\left(\frac{\zeta\omega'(\eta\zeta)\overline{\omega(\eta\zeta)}}{1-|\omega(\eta\zeta)|^{2}}\right),
$$
which shows that
\beq\label{eq-11}
\mbox{Re}\left(\frac{\zeta h''(\eta\zeta)}{h'(\eta\zeta)}\right)&\leq&
\mbox{Re}\left(\frac{\zeta\omega'(\eta\zeta)\overline{\omega(\eta\zeta)}}{1-|\omega(\eta\zeta)|^{2}}\right)+\frac{\nu}{1-\eta^{2}}\\
\nonumber
&\leq&\frac{|\omega'(\eta\zeta)||\overline{\omega(\eta\zeta)|}}{1-|\omega(\eta\zeta)|^{2}}+\frac{\nu}{1-\eta^{2}},
\eeq
where $\zeta\in\partial\mathbb{D}.$ By Schwarz-Pick's lemma, we obtain
\be\label{eq-12}
|\omega'(\eta\zeta)|\leq\frac{1-|\omega(\eta\zeta)|^{2}}{1-\eta^{2}}.
\ee
By (\ref{eq-11}) and (\ref{eq-12}), we have
$$\mbox{Re}\left(\frac{\zeta h''(\eta\zeta)}{h'(\eta\zeta)}\right)\leq\frac{1+\nu}{1-\eta^{2}}.
$$
Choosing $\lambda\in(0,1-\nu)$, there is an $r_{1}\in[r_{0},1)$ such that
\be\label{eq-y1}
\mbox{Re}\left(\frac{\zeta
h''(\eta\zeta)}{h'(\eta\zeta)}\right)<\frac{2\eta-2\lambda}{1-\eta^{2}}
~\mbox{ for all $\zeta\in\partial\mathbb{D}$},
\ee
when $\eta\in[r_{1},1)$. For $0\leq r_{1}\leq r\leq\rho<1$, by (\ref{eq-y1}), we find that
\beq \nonumber
\log\left[\frac{(1-\rho^{2})|h'(\rho\zeta)|}{(1-r^{2})|h'(r\zeta)|}\right]&=&\int_{r}^{\rho}\bigg[\mbox{Re}\left(\frac{\zeta
h''(\eta\zeta)}{h'(\eta\zeta)}\right)-\frac{2\eta}{1-\eta^{2}}\bigg]d\eta\\
\nonumber
&<&-2\lambda\int_{r}^{\rho}\frac{d\eta}{1-\eta^{2}}\\
\nonumber
&=&-\lambda\log\left(\frac{1+\rho}{1+r}\cdot\frac{1-r}{1-\rho}\right),
\eeq
which implies that
\be\label{eq-y2}
\left|\frac{h'(\rho\zeta)}{h'(r\zeta)}\right|<\left(\frac{1+r}{1+\rho}\right)^{1+\lambda}
\left(\frac{1-\rho}{1-r}\right)^{\lambda-1}\leq\left(\frac{1-\rho}{1-r}\right)^{\lambda-1}.
\ee
By (\ref{eq-y2}), we get
\beq\label{eq-y3}
\|D_{f}(\rho\zeta)\|&\leq&\frac{2K}{1+K}|h'(\rho\zeta)|<\frac{2K}{1+K}
|h'(r\zeta)|\left(\frac{1-\rho}{1-r}\right)^{\lambda-1}\\ \nonumber
&\leq&\frac{2K}{1+K} \|D_{f}(r\zeta)\|\left(\frac{1-\rho}{1-r}\right)^{\lambda-1}.
\eeq
In order to apply Theorem \ref{cor-1} and then to conclude that
$\Omega_{1}=f(\mathbb{D})$ is a radial John disk, we will use some
proof techniques as in the proof of \cite[Theorem 3.7]{KH} to remove
the restriction $r\geq r_{1}$ above. For $0\leq r_{1}\leq
r\leq\rho<1$, by (\ref{eq-2.0}) and (\ref{eq-y3}), we see that there
is a constant $c(\lambda)>1$ such that
\be\label{eq-u1}
\sigma_{\ell}(w)\leq c(\lambda)d_{\Omega_{1}}(w),
\ee
where $w_{1}=f(\rho\zeta)$, $w=f(r\zeta)$ and $\gamma=f([0,\rho\zeta])$.
It follows from (\ref{eq-u1}) that
\be\label{eq-u2}
\diam(\gamma[w_{1},w])\leq c(\lambda)d_{\Omega_{1}}(w).
\ee

Now we consider the case: $0\leq r\leq r_{1}\leq\rho<1$. Let
$\delta_{0}=\dist\big(f(\overline{\mathbb{D}}_{r_{1}}),\partial\Omega_{1}\big)$
denote  the Euclidean distance from $f(\overline{\mathbb{D}}_{r_{1}})$ to the boundary
$\partial\Omega_{1}$ of $\Omega_{1}$ and let
$\lambda_{0}=\diam\big(f(\overline{\mathbb{D}}_{r_{1}})\big)$. Then
\be\label{eq-u3}
\delta_{0}>0~\mbox{and}~\lambda_{0}<\infty.
\ee
For $0\leq r\leq r_{1}\leq\rho<1$, by the triangle inequality, (\ref{eq-u2}) and (\ref{eq-u3}), we get
\beq \nonumber
\diam(\gamma[w,w_{1}])&\leq&\diam(\gamma[w,w_{0}])+\diam(\gamma[w_{0},w_{1}])\\
\nonumber &\leq&\lambda_{0}+c(\lambda)d_{\Omega_{1}}(w_{0})\\
\nonumber&\leq&\lambda_{0}+c(\lambda)(\lambda_{0}+\delta_{0})\\
\nonumber &\leq&(c(\lambda)+c')\delta_{0}\\
\nonumber &\leq&(c(\lambda)+c')d_{\Omega_{1}}(w),
\eeq
where $w_{1}=f(\rho\zeta)$, $w=f(r\zeta)$, $w_{0}=f(r_{1}\zeta)$ and $c'=(1+c(\lambda))\lambda_{0}/\delta_{0}$.

The remaining case when $0\leq r\leq\rho \leq r_{1}<1$ is treated similarly. Therefore, for $0\leq r\leq\rho<1$,
there is a constant $c_{2}>1$ such that
$$\diam\big(\gamma[w,w_{1}]\big)\leq c_{2}d_{\Omega_{1}}(w),
$$
which implies that $\mbox{car}_{d}(\gamma,c_{2})\subset\Omega_{1}$ (cf. \cite{NV}), where
$$\mbox{car}_{d}(\gamma,c_{2})=\bigcup\Big\{\mathbb{D}\big(w,\diam\big(\gamma[w,w_{1}]\big)/c_{2}\big):~w\in\gamma\setminus\{f(0),w_{1}\}\Big\}.
$$
It follows from   \cite[Theorem 2.16]{NV}  and \cite[Part 2.26 in P.17]{NV} that $\Omega_{1}$ is a John disk.  For the definition of the diameter of
$c$-carrot, denoted by $\mbox{car}_{d}(\gamma, c)$, we refer to \cite{NV}. The proof of the theorem is complete.
 \hfill $\Box$

\bigskip

{\bf Acknowledgements:}  This research was partly supported by the
National Natural Science Foundation of China (No. 11401184 and No.
11571216), the Hunan Province Natural Science Foundation of China
(No. 2015JJ3025), the Excellent Doctoral Dissertation of Special
Foundation of Hunan Province (higher education 2050205), the
Construct Program of the Key Discipline in Hunan Province (No.
[2011] 76). The second author is currently on leave from
  Indian Institute of Technology Madras.


\normalsize


\begin{thebibliography}{99}

\bibitem{AW} {\sc L. V. Ahlfors and G. Weill,}
A uniqueness theorem for Beltrami equstions,
\textit{Proc. Amer. Math. Soc.,} {\bf 13} (1962), 975--978.

\bibitem{BP} {\sc J. Becker and Ch. Pommerenke,}
Schlichtheitskriterien und jordangebiete,
\textit{J. Reine Angew. Math.,} {\bf 354} (1984), 74--94.

\bibitem{CPRW} {\sc Sh. Chen, S. Ponnusamy, A. Rasila and X. Wang,}
Linear connectivity, Schwarz-Pick lemma and univalency criteria for planar harmonic mappings,
\textit{Acta Math. Sinica}, 13 pages, To appear.\\
{\tt http://arxiv.org/pdf/1404.4155v1.pdf}

\bibitem{CDO} {\sc M. Chuaqui, P. Duren and B. Osgood,}
The Schwarzian derivative for harmonic mappings,
\textit{J. Anal. Math.,} {\bf 91} (2003), 329--351.

\bibitem{COC} {\sc M. Chuaqui, B. Osgood and Ch. Pommerenke,}
John domains, quasidisks and the Nehari class,
\textit{J. Reine Angew. Math.,} {\bf 471} (1996), 77--114.

\bibitem{Clunie-Small-84} {\sc J. G.  Clunie and T. Sheil-Small,}
Harmonic univalent functions,
\textit{Ann. Acad. Sci. Fenn. Ser. A I Math.,} {\bf 9} (1984), 3--25.

\bibitem{Co}{\sc F. Colonna,}
The Bloch constant of bounded harmonic mappings, \textit{Indiana
Univ. Math. J.}, {\bf 38}(1989), 829--840.

\bibitem{Du} {\sc P. Duren,}
{\it Harmonic mappings in the plane,} Cambridge Univ. Press, 2004.

\bibitem{KH} {\sc K. Hag and P. Hag,}
John disks and the pre-Schwarzian derivative,
\textit{Ann. Acad. Sci. Fenn.  Math.,} {\bf 26} (2001), 205--224.

\bibitem{HM} {\sc R. Hern\'andez and M. J. Mart\'in,}
Pre-Schwarzian and Schwarzian derivatives of harmonic mappings,
\textit{J. Geom. Anal.,} {\bf 25} (2015), 64--91.

\bibitem{John} {\sc F. John,}
Rotation and strain,
\textit{Comm. Pure Appl. Math.,} {\bf 14} (1961), 391--413.

\bibitem{K} {\sc D. Kalaj,}
Muckenhoupt weights and Lindel\"of theorem for harmonic mappings,
\textit{Adv. Math.,} {\bf 280} (2015), 301--321.


\bibitem{MM} {\sc M. Kne{\rm $\breve{z}$}evi\'c and M. Mateljevi\'c,}
On the quasi-isometries of harmonic quasiconformal mappings,
\textit{J. Math. Anal. Appl.,} {\bf 334} (2007), 404--413.

\bibitem{LV} {\sc O. Lehto and K. I. Virtanen,}
{\it Quasiconformal mappings in the plane}, Springer Verlag, 1973.

\bibitem{Lewy} {\sc H. Lewy,}
On the non-vanishing of the Jacobian in certain one-to-one mappings,
\textit{Bull. Amer. Math. Soc.,} {\bf 42} (1936), 689--692.

\bibitem{Mi} {\sc M. Mateljevi\'c,}
Distrotion of quasiregular mappings and equivalent norms on Lipschitz-type spaces,
\textit{Abstr. Appl. Anal.,} Volume 2014, Article ID 895074, 20 pages. 

\bibitem{NV} {\sc R. N\"akki and J. V\"ais\"al\"a,}
John disks,
\textit{Exposition Math.,} {\bf 9} (1991), 3--43.

\bibitem{Po} {\sc Ch. Pommerenke,}
One-sided smoothness conditions and conformal mapping,
\textit{J. London Math. Soc.,} {\bf 26} (1982), 77--88.

\bibitem{Po1} {\sc Ch. Pommerenke,}
{\it Boundary behaviour of conformal maps}, Springer-Verlag, 1992.

\bibitem{SaRa2013}    {\sc S. Ponnusamy and A. Rasila},
Planar harmonic and quasiregular mappings.
Topics in Modern  Function Theory (Editors. St. Ruscheweyh and S. Ponnusamy): Chapter in
CMFT, RMS-Lecture Notes Series No. 19, 2013, pp. 267--333. 

\bibitem{Small} {\sc  T. Sheil-Small,}
Constants for planar harmonic mappings, \textit{J. London Math.
Soc.,} {\bf 42} (1990), 237--248.

\bibitem{Va} {\sc J. V\"ais\"al\"a,}
{\it Lectures on $n$-dimensional quasiconformal mappings,}
Springer-Verlag, Berlin, Heidelberg, New York,  xiv, 144pp, 1971.

\bibitem{V} {\sc M. Vuorinen,}
{\it Conformal geometry and quasiregular mappings,}
\textit{Lecture Notes in Math.} Vol. 1319, Springer-Verlag, 1988.

\end{thebibliography}
\end{document}